\documentclass[12pt]{article}
\usepackage{latexsym,amsmath,amssymb,amsfonts,amstext,amsthm}

 \numberwithin{equation}{section}

\begin{document}
\title{\bf  Mysterious Properties of the Point at Infinity
}

\author{ Saburou Saitoh \\ Institute of Reproducing Kernels
}


\maketitle
 {\bf Abstract:}  From the viewpoints of the division by zero $1/0=0/0=z/0=0$ and the division by zero calculus, we will examine the mysterious properties of the point at infinity in the sense of the Alexandroff one compactification of  the complex plane which is realized by the stereographic projection.


\section{Introduction}

 The purposes of this paper are to introduce the very fundamental new concept of the point at infinity and to propose the related open problems.

 The division by zero $1/0=0/0=z/0$ itself will be quite clear and trivial with several {\bf  natural extensions} of the fractions against the mysteriously long history (\cite{romig}), as we can see from the concepts of the Moore-Penrose generalized inverse or the Tikhonov regularization method to the fundamental equation $ az =b$ whoes solutions lead to the definition $z= b/a$.
  
  However, the result will show that
      for the elementary mapping
\begin{equation}
W = \frac{1}{z},
\end{equation}
the image of $z=0$ is $W=0$ ({\bf should be defined from the form}). This fact seems to be a curious one in connection with our well-established popular image for the  point at infinity on the Riemann sphere. As the representation of the point at infinity of the Riemann sphere by the
zero $z =  0$, we will see some delicate relations between $0$ and $\infty$ which show a strong
discontinuity at the point of infinity on the Riemann sphere. We did not consider any value of the elementary function $W =1/ z $ at the origin $z = 0$, because we did not consider the division by zero
$1/ 0$ in a good way. Many and many people consider its value by the limiting like $+\infty $ and  $- \infty$ or the
point at infinity as $\infty$. However, their basic idea comes from {\bf continuity} with the common and natural sense or
based on the basic idea of Aristotle.  --
  For the related Greece philosophy, see \cite{a,b,c}. However, as the division by zero we will consider the value of
the function $W =1 /z$ as zero at $z = 0$. We will see that this new definition is valid widely in
mathematics and mathematical sciences, see  (\cite{mos,osm}) and the cited references for example. Therefore, the division by zero will give great impacts to calculus, Euclidean geometry, differential equations, analytic geometry, complex analysis and physics  in the undergraduate level and to our basic ideas for the space and universe. Here, we would like to refer to some mysterious properties of the point at infinity.

\section{Logical background}
As stated in the introduction, our results are different from the very classical results and ideas since  Aristotle and
 Eulclid. Therefore, we will check and  confirm simply our new mathematics.

At first, the number system containing the division  by zero is established as the {\bf Yamada field} in \cite{msy} by a minor arrangement of the well-established complex number field. The Yamada field is a natural and simple extension  of the complex field containing the division by zero.

The {\bf uniqueness theorem} for the generalized fraction
(division) is established by Takahasi \cite{kmsy} under the very general  assumption of the product $(a/b) \cdot (c/d) = (ac)/(bd)$. If  this property is not satisfield, then we will not be able to find any fundamental meanings of the fractions.

Our division  by zero is also given as the {\bf Moore-Penrose generalized inverse} of the fundamental equation $az=b$ whose solution is represented by the generalized fraction $z=b/a$. The  Moore-Penrose generalized inverse is a well-established general concept.

We gave many and many clear interpretations  and applications of our division by zero in the cited references.

We will be able to confirm our division by zero as the natural mathematics as the number system.

However, we will need the concept of the {\bf division by zero calculus} for applying the division by zero to {\bf functions}. This will be given simply in the following way.

For any  Laurent expansion around $z=a$,

\begin{equation}\label{eq2}
f(z) = \sum_{n=-\infty}^{- 1} C_n (z - a)^n +  C_0 +  \sum_{n= 1}^{\infty} C_n (z - a)^n,
\end{equation}
we obtain the identity, by the division by zero

\begin{equation}\label{eq3}
f(a) =  C_0. 
\end{equation}

\medskip

For the correspondence \eqref{eq3} for the function $f(z)$, we will call it 
{\bf the division by zero calculus}. By considering the derivatives in 
\eqref{eq2}, we can define any order derivatives of the function $f$ at the 
singular point $a$.

For functions, we can, in general, consider the applications of the division by zero in many ways, however, the above division by zero calculus is applicable in many cases, however, for the results obtained {\bf we have to check} their senses. Indeed, we have the cases that we have good results and nonsense results. 

For example, for the simple example for the line equation on the $x, y$ plane
$$
 ax + by + c=0
$$
we have, formally
$$
x + \frac{by + c}{a} =0,
$$
and so, by the division by zero, we have, for $a=0$, the reasonable result
$$
x = 0.
$$

Indeed, for the equation $y =mx$, from
$$
\frac{y}{m} =x,
$$
we have, by the division by zero, $x=0$ for $m=0$. This gives the case $m=\infty$ of the gradient of the line.
-- This will mean that the equation $y =mx$ represents the general line through the origin in this sense. -- This method was applied in many cases, for example see \cite{os, o}.

However, from
$$
\frac{ax + by}{c} + 1 =0,
$$
for $c=0$, we have the contradiction, by the division by zero
$$
1 =0.
$$

Meanwhile, note that for the function $f(z) = z + \frac{1}{z}$, $f(0) =0$, however, for the function
$$
f(z)^2 = z^2 + 2 + \frac{1}{z^2},
$$
we have $f^2 (0) = 2$. Of course,
$$
f(0) \cdot f(0) =\{f(0)\}^2 = 0.
$$

We can consider many ways
applying the division by zero  to functions. For many concrete examples,  see  (\cite{mos,osm}) and the papers cited in the references. In particular, note that we can consider the division by zero for more general functions that are not restricted to analytic functions.

Before ending the background for this paper, we have to refer to the very important facts:
 
We have even the formal contradiction for the very classical result that the point at infinity is represented by $\infty$ and it is represented by zero.
We were able to establish the fundamental relation between the point at infinity and $0$,  see \cite{mos}.

The inversion with respect to a circle of the center of the circle is given by the center of the circle, not the point at infinity. A line  may be looked as a circle with radius zero and with center at the origin. See
(\cite{mos,osm}).

\section{Many points at infinity?}

When we consider a circle with center $P$, by the inversion with respect to the circle,  the points of a neighborhood at the point $P$ are mapped to a neighborhood around the point at infinity except the point $P$.  This property is independent of the radius of the circle. It looks that the point at infinity is depending on the center $P$. This will mean that there exist many points at infinity,  in a sense.

\section{Stereographic projection}

The point at infinity may be realized by the stereographic projection as well known. However, the projection is dependent on the position of the sphere (the plane coordinates).  Does this mean that there exist many points at infinity?

\section{Laurent expansion}

From the definition of the division  by zero calculus,
 we see that if there exists a negative $n$ term in (2.1)
$$
\lim_{z\to a} f(z) = \infty,
$$
however, we have (2.2). The values at the point $a$ have many values, that are all complex numbers. At least, {\bf in this sense}, we see that we have many points as the point of infinity.


In the sequel, we will show typical points at infinity.

\section{Diocles' curve of Carystus (BC 240? - BC 180?)}

The beautiful curve

$$ y^2 = \frac{x^3}{2a - x}, \quad a > 0
$$
 is considered by Diocles. By setting $ X = \sqrt{2a -x}$ we have
$$
y = \pm \frac{x^{(3/2)}}{\sqrt{2a -x}} =\pm \frac{(2a -X^2)^{(3/2)}}{X}.
$$
Then, by the division by zero calculus at $X=0$, we have a reasonable value $0$.

Meanwhile, for the function $\frac{x^3}{2a - x}$, we have $-12a^2$, by the division by zero calculus at $x = 2a$. This leads to a wrong value.

\section{Nicomedes' curve (BC 280 - BC 210)}

The very interesting curve
$$
r= a + \frac{b}{\cos \theta}
$$
is considered by Nicomedes from the viewpoint of the $1/3$ division of an angle. That has  very interesting geometrical meanings.
For the case $\theta = \pm (\pi/2)$, we have $r = a$, by the division by zero calculus.

Of course, the function is symmetric for $\theta = 0$, however, we have a mysterious value $r=a$, for $\theta = \pm (\pi/2)$. Look the beautiful graph of the function.

\section{Newton's curve (1642 - 1727)}

Meanwhile, for the famous Newton curve
$$
y = a x^2 + bx + c + \frac{d}{x} \quad (a, d \ne 0),
$$
of course, we have $y(0) =c$.

Meanwhile, in the division by zero calculus, the value  is determined by the information around any analytical point for an analytic function, as we see from the basic property of analytic functions.

At this moment, the properties of the values of  analytic functions at isolated singular points are mysterious, in particular, in the  geometrical sense.

\section{Basic meanings of  values at isolated singular points of analytic functions}
 
 Since the values of analytic functions at isolated singular points were given by the coefficients $C_0$ of the Laurent expansions (2.1) as the division by zero calculus. Therefore,  their values may be considered as arbitrary ones by any sift of the image complex plane.  Therefore, we can consider the  values as zero in any Laurent expansions by  shifts, as normalizations. However, if  the Laurent expansions are determined by another normalizations, then the  values will have their senses. We will examine such properties for the Riemann mapping function.  
 
 Let $D$ be a simply-connected domain containing the point at infinity having at least two boundary points. Then, by the celebrated theorem of Riemann, there exists a uniquely determined conformal mapping with a series expansion
 \begin{equation}
 W = f(z) = C_1 z + C_0 + \frac{C_{-1}}{z} + \frac{C_{-2}}{z^2} +\dots,   C_1 > 0,
 \end{equation}
 at the point at infinity which maps the domain $D$ onto the exterior $|w| >1$ of the unit disc on the complex $W$ plane.  We can normalize (9.1) as follows:
 \begin{equation}
 \frac{f(z) }{C_1}  =  z +  \frac{C_0}{C_1} + \frac{C_{-1}}{C_1z} + \frac{C_{-2}}{C_1z^2} +\dots.
 \end{equation}
 Then, this function $\frac{f(z) }{C_1}$ maps $D$ onto the exterior of a circle of radius 
 $1/C_1$ and so, it is called the {\bf mapping radius} of $D$ (\cite{bs}, 80 page; \cite{tsuji}, 83 page). 

Meanwhile, from the normalization
\begin{equation}
 f(z) - C_0 =  C_1 z +    \frac{C_{-1}}{z} + \frac{C_{-2}}{z^2} +\dots,
 \end{equation}
 by the natural shift  $C_0 $ of the image plane, the  unit circle is mapped to the unit  circle with center $C_0 $. Therefore, $C_0 $ may be called as {\bf mapping center} of  $D$.
 The function $f(z)$ takes the value $C_0 $ at the point at infinity in the sense of the division by zero calculus and now we have its natural sense as the mapping center of $D$. We have  considered the value of the function $f(z)$ as infinity at the point at infinity, however, practically it was the value $C_0 $. This will mean that in a sense the value $C_0 $ is the farthest point from the point at infinity or from the image domain with the strong discontinuity. -- Recall the mapping property of the fundamental function $W=1/z$ of the unit disc $|z| <1$.
 
 The properties of mapping radius were investigated deeply in conformal mapping theory like estimations, extremal properties and meanings of the values, however, it seems that there is no information on the property of mapping center. See many books on conformal mapping theory or analytic function theory. See \cite{tsuji} for example.

\section{Unbounded, however, bounded}

We will consider the high
$$ 
y = \tan \theta,      \quad  0 \le \theta \le \frac{\pi}{2}
$$
on the line $x=1$.
Then, the high $y$ is unbounded, however, the high line (gradient) can not be extended beyond the $y$ axis.  The restriction is given by $0 = \tan (\pi/2)$.

Recall the stereographic projection of the complex plane. The points on the plane can be expanded in an unbounded way, however, all the points on the complex plane have to be corresponded to the points of the Riemann sphere. The restriction is the point at infinity which corresponds to the north pole of the Riemann sphere and the point at infinity is represented by 0.

\bibliographystyle{plain}

\noindent
Saburou Saitoh \\ 
Institute of Reproducing Kernels\\
Kawauchi-cho, 5-1648-16, Kiryu 376-0041, Japan\\
{\it  E-mail address}: kbdmm360@yahoo.com.jp

\end{document}